# Variance and the Inequality of Arithmetic and Geometric Means

Burt Rodin

**Abstract**

There have been a number of recent papers devoted to generalizations of the classical AM-GM inequality. Those generalizations which incorporate *variance* have been the most useful in applications to economics and finance. In this paper we prove an inequality which yields best possible upper and lower bounds for the geometric mean of a sequence solely in terms of its arithmetic mean and its variance. A particular consequence is the following: among all positive sequences having given length, arithmetic mean and nonzero variance, the geometric mean is maximal when all terms in the sequence except one are equal to each other and are less than the arithmetic mean.

**Introduction.** Roughly speaking, the discrepancy between the arithmetic and geometric means of a finite sequence tends to increase as the sequence deviates more and more from being constant. The literature contains several generalizations of the classical arithmetic-geometric mean inequality; they differ, in part, by using different measures for the deviation of the sequences from constancy. *Variance* (or *standard deviation* ) is a mathematically natural measure of the deviation of a sequence from constancy. In addition, as noted in Aldaz [1,2], Becker [4], Estrada [6], and Markowitz [9], variance is the most useful such measure from the point of view of economics and finance (as Markowitz [9] points out, investors are made aware of the arithmetic mean and variance of a portfolio, but there is a need for them to estimate the geometric mean since that is the portfolio's likely long term return; cf. Remark 2 below). Theorem 1 in this paper gives bounds for the geometric mean depending solely on the arithmetic mean and variance: these mean-variance bounds are best possible. A discussion of related previous results is given in Section 2. Corollary 1 in Section 2 yields an upper bound, depending only on variance, for the numerical difference between the arithmetic and geometric means (cf. Aldaz [2]).

## Section 1

Let $x_1, x_2, \ldots , x_n$ be a sequence of n real numbers. The (arithmetic) mean $\mu$ and variance $\sigma^2$ are defined as follows :

(1)
$$\mu = \frac{1}{n}\sum_{i=1}^{n} x_i, \quad \sigma^2 = \frac{1}{n}\sum_{i=1}^{n}(x_i - \mu)^2.$$

This notation always implies that σ is the standard deviation, i.e., the nonnegative square root of the variance. Thus *mean* will always refer to the arithmetic mean; the quantity $(x_1 x_2 \cdots x_n)^{1/n}$ will be referred to by its complete name, *geometric mean*.

It is clear that, for fixed mean $\mu > 0$, if $\sigma^2$ is sufficiently small then each $x_i$ will necessarily be positive. Also, for fixed mean $\mu > 0$, if all the $x_i$ are positive then the variance of the sequence cannot be too large. The precise conditions for the mean and variance with these two properties are given in the following lemma.

**Lemma 1.** *Let $x_1, x_2, \ldots, x_n$ be a sequence of $n \geq 2$ real numbers with mean $\mu > 0$ and variance $\sigma^2$. (a) If $\sigma/\mu < 1/\sqrt{n-1}$ then all terms of the sequence are necessarily positive. (b) If all terms of the sequence are positive then $\sigma/\mu < \sqrt{n-1}$.*

**Proof.** Let $\mu > 0$ be fixed. Let S be the (n-1)-simplex

(2)    $S = \{ (x_1, x_2, \ldots, x_n) \in \mathbf{R}^n : x_1 + x_2 + \ldots + x_n = n\mu$ and $x_i \geq 0$ for $i=1,\ldots,n\}$.

The variance $\sigma^2$ of the coordinates of a point $(x_1, x_2, \ldots, x_n)$ of S is related to the distance r from that point to the centroid $C_0 = (\mu, \mu, \ldots, \mu)$ of *S* by $r^2 = n\sigma^2$.

Let $r_1$ be the distance from $C_0$ to a nearest boundary point of S and let $r_2$ be the distance from $C_0$ to a furthest boundary point of S. If $(x_1, x_2, \ldots, x_n)$ has mean $\mu$ and if the distance from $(x_1, x_2, \ldots, x_n)$ to $C_0$ is $\leq r_1$, then each coordinate $x_i$ must be nonnegative. Similarly, if $x_1, x_2, \ldots, x_n$ has mean $\mu$ and if each $x_i$ is nonnegative, then the distance from $(x_1, x_2, \ldots, x_n)$ to $C_0$ must be $\leq r_2$.

The boundary points of S nearest to $C_0$ are the centroids of each (n-2)-face of S, for example, the point $(0, n\mu/(n-1), \ldots, n\mu/(n-1))$. The distance $r_1$ from $C_0$ to such a nearest boundary point satisfies $r_1^2 = \mu^2 n/(n-1)$. Therefore, if a sequence of $x_1, x_2, \ldots, x_n$ with mean $\mu$ has variance $\sigma^2 < r_1^2/n = \mu^2/(n-1)$ then all terms of that sequence are necessarily positive. This proves (a) of Lemma 1.

Similarly, the boundary points of S furthest from $C_0$ are the vertices of S. The distance $r_2$ from $C_0$ to a vertex of S satisfies $r_2^2 = \mu^2 n(n-1)$. Therefore, if a sequence $x_1, x_2, \ldots, x_n$ with mean μ and variance σ² has all term positive then σ² < $r_2^2$/n = μ²(n-1). This proves (b) of Lemma 1.

**Theorem 1.** *Let* n ≧ 2. *Let* $x_1, x_2, \ldots, x_n$ *be real numbers with mean* μ > 0 *and variance* σ².

(a) *If*  $0 \le \sigma/\mu < 1/\sqrt{n-1}$  *then each* $x_i$ *is positive and*

(3) $\quad (\mu - \sigma\sqrt{n-1})\left(\mu + \dfrac{\sigma}{\sqrt{n-1}}\right)^{n-1} \le x_1 x_2 \cdots x_n \le (\mu + \sigma\sqrt{n-1})\left(\mu - \dfrac{\sigma}{\sqrt{n-1}}\right)^{n-1}.$

*The upper and lower bounds in (3) are sharp.*

(b) *If every term of the sequence* $x_1, x_2, \ldots, x_n$ *is positive then* $0 \le \sigma/\mu < \sqrt{n-1}$ *and the inequalities (3) continue to hold. The upper bound is again sharp. In the subrange* $1/\sqrt{n-1} < \sigma/\mu < \sqrt{n-1}$ *the lower bound expression in (3) becomes negative and should be replaced by 0; with that understanding the lower inequality will then be best possible for the entire range* $0 \le \sigma/\mu < \sqrt{n-1}$.

**Remark 1.** Up to a change in the order of the terms, the sequences which make (3) an equality are the following. For  $0 \le \sigma/\mu < \sqrt{n-1}$  the upper bound is attained when

(4) $\quad\quad\quad\quad x_1 = x_2 = \ldots = x_{n-1} = \mu - \dfrac{\sigma}{\sqrt{n-1}}$  and  $x_n = \mu + \sigma\sqrt{n-1}$.

For  $0 \le \sigma/\mu < 1/\sqrt{n-1}$  the lower bound is attained when

(5) $\quad\quad\quad\quad x_1 = x_2 = \ldots = x_{n-1} = \mu + \dfrac{\sigma}{\sqrt{n-1}}$  and  $x_n = \mu - \sigma\sqrt{n-1}$ ;

for  $1/\sqrt{n-1} \le \sigma/\mu < \sqrt{n-1}$  there is no minimum among positive sequences with the given μ,σ but the infimum is 0.

**Proof of Theorem 1.** Let n ≧ 2, μ > 0, and σ² be fixed. Let **x** = $(x_1, x_2, \ldots, x_n)$ be a real n-vector and let

(6) $\quad G(\mathbf{x}) = x_1 x_2 \cdots x_n, \quad A(\mathbf{x}) = x_1 + x_2 + \cdots + x_n, \quad V(\mathbf{x}) = (x_1-\mu)^2 + (x_2-\mu)^2 + \ldots + (x_n-\mu)^2.$

Consider the problem: maximize or minimize G(**x**) subject to the constraints

(7) $\qquad\qquad A(\mathbf{x})= n\mu$ and $V(\mathbf{x})=n\sigma^2$.

We shall refer to this as the *max-min problem*. By a *critical point* for this problem we mean a point **x** which satisfies constraints (7) and where **grad** $G(\mathbf{x})$ is in the space spanned by **grad** $A(\mathbf{x})$ and **grad** $V(\mathbf{x})$. The method of Lagrange multipliers asserts that the solutions to the max-min problem—which exist by compactness—will be found among the values of G at the critical points. First we find all these critical points. Then we will consider the restrictions in the theorem regarding positivity and bounds on $\sigma/\mu$.

Thus **x** will be a critical point if equations (7) are satisfied and if there exist numbers $\lambda_1, \lambda_2$ such that

(8) $\qquad\qquad$ **grad** $G(\mathbf{x}) = \lambda_1$ **grad** $A(\mathbf{x}) + \lambda_2$ **grad** $V(\mathbf{x})$.

We have

$\qquad\qquad$ **grad** $G(\mathbf{x}) = (x_2 x_3 \cdots x_n, \ldots, x_1 x_2 \cdots x_{n-1})$,

(9) $\qquad\qquad$ **grad** $A(\mathbf{x}) = (1,1,\ldots,1)$,

$\qquad\qquad$ **grad** $V(\mathbf{x}) = 2(x_1-\mu, x_2-\mu, \ldots, x_n-\mu)$.

If we multiply each side of (8) by $x_i$ and then equate the i-th components on each side we obtain via the three equations (9) the n scalar equations

(10) $\qquad\qquad x_1 x_2 \cdots x_n = (\lambda_1 - 2\mu\lambda_2)x_i + 2\lambda_2 x_i^2, \qquad (i=1,2,\ldots,n)$.

If **grad** $G(\mathbf{x})=\mathbf{0}$ then at least two of the coordinates of **x** are 0; the converse is also true. Suppose **grad** $G(\mathbf{x}) \neq \mathbf{0}$. Then $\lambda_1, \lambda_2$ are not both zero. Equations (10) then show that all of the n ordered pairs $(x_i, x_i^2)$ lie on a line

(11) $\qquad\qquad (\lambda_1 - 2\mu\lambda_2)\, x + 2\lambda_2\, y = \text{constant};$

of course they also lie on the parabola $y=x^2$. Therefore there are at most two distinct values in the set $\{x_1, x_2, \ldots, x_n\}$.

We have seen that a point $\mathbf{x} = (x_1, x_2, \ldots, x_n)$ is a critical point for the max-min problem only if **x** has at most two distinct coordinates or else **grad** $G(\mathbf{x})=\mathbf{0}$. The case of one distinct coordinate $x_1 = x_2 = \cdots = x_n$ occurs if and only if $\sigma=0$. In this case the inequalities (3) become trivial equalities.

Consider a critical point **x** such that **grad** $G(\mathbf{x}) \neq \mathbf{0}$ and such that the coordinates of **x** have exactly two distinct values; denote these two values by $a$ and $b$, with $b < \mu < a$.

Suppose the value *a* occurs *i* times and *b* occurs *j* times, where $0 < i,j < n$ and $i+j=n$. The constraints (7) require that

(12) $\qquad ia+jb=n\mu$ and $i(a-\mu)^2 +j(b-\mu)^2 =n\sigma^2$.

To express *a* and *b* in terms of $i,j,n,\mu,\sigma$, solve the first equation in (12) for *b*, substitute the solution into the second equation, and obtain

(13) $$i(a-\mu)^2 = j\sigma^2.$$

Since $a>\mu$,

(14) $$a=\mu+\sigma\sqrt{\frac{j}{i}}.$$

Now the first equation in (12) yields

(15) $$b=\mu-\sigma\sqrt{\frac{i}{j}}.$$

Given *i*, $1 \leq i \leq n-1$, there are $\binom{n}{i}$ points which have *i* coordinates *a* given by (14) and $j=n-i$ coordinates *b* given by (15); these will be called *critical points of type i*. If $x_i$ is a critical point of type *i* then the corresponding critical value of G is

(16) $$G(x_i) = \left(\mu+\sigma\sqrt{\frac{j}{i}}\right)^i \cdot \left(\mu-\sigma\sqrt{\frac{i}{j}}\right)^j, \qquad (i=1,\ldots,n-1;\ j=n-i).$$

We now want to order the n-1 critical values in (16) according to magnitude.

Let $t=\sigma/\mu$. Define $P_i(t)$ by

(17) $$P_i(t) = \frac{G(x_i)}{\mu^n}.$$

Each $P_i(t)$ ($1 \leq i \leq n-1$) can be considered as a polynomial in *t* of degree n:

(18) $$P_i(t) = \left(1+t\sqrt{\frac{j}{i}}\right)^i \cdot \left(1-t\sqrt{\frac{i}{j}}\right)^j, \qquad (i=1,\ldots,n-1;\ j=n-i).$$

**Lemma 2**. *Let i, j, and $P_i(t)$ be given by Equation (18). If $1 \leq i \leq n-2$ then*

(19) $$P_i(t) > P_{i+1}(t) \text{ for } 0 < t < \sqrt{\frac{j-1}{i+1}}.$$

*Consequently,*

(20) $$P_{n-1}(t) < \cdots < P_2(t) < P_1(t) \text{ for } 0 < t < \frac{1}{\sqrt{n-1}}.$$

**Proof of Lemma 2.** From equation (18) we find that for $1 \leq i \leq n-1$

(21) $$\frac{d}{dt} \log P_i(t) = \frac{-nt}{(1+t\sqrt{j/i})(1-t\sqrt{i/j})}.$$

Thus $P_i(t)$ decreases from 1 to 0 as t goes from 0 to $\sqrt{j/i}$. We have

(22) $$\log P_i(t) = \int_0^t \frac{-n\tau \, d\tau}{(1+\tau\sqrt{j/i})(1-\tau\sqrt{i/j})}, \quad (0 < t < \sqrt{j/i}).$$

According to the representation (22), in order to prove that for $1 \leq i \leq n-2$

(23) $$\log P_i(t) > \log P_{i+1}(t) \text{ for } 0 < t < \sqrt{\frac{j-1}{i+1}}$$

it suffices to show that for $0 < \tau < \sqrt{(j-1)/(i+1)}$

(24) $$\frac{-n\tau}{(1+\tau\sqrt{j/i})(1-\tau\sqrt{i/j})} > \frac{-n\tau}{(1+\tau\sqrt{(j-1)/(i+1)})(1-\tau\sqrt{(i+1)/(j-1)})}.$$

For τ in this range the factors in the denominators of (24) are positive and the inequality (24) can be algebraically simplified to become

(25) $$\sqrt{\frac{j-1}{i+1}} - \sqrt{\frac{i+1}{j-1}} < \sqrt{\frac{j}{i}} - \sqrt{\frac{i}{j}}.$$

Replace *j* by n-*i* in (25); the result can be written as

(26) $$\frac{n-2i}{\sqrt{i(n-i)}} > \frac{n-2(i+1)}{\sqrt{(i+1)(n-(i+1))}}.$$

It is evident that the left hand side of (26) is a strictly decreasing function of a real variable $i$ in the interval $0<i<n$ since its derivative with respect to $i$ is negative. Therefore (26) is valid for integers $i$ in the range $1 \leq i \leq n-2$, and (23) follows. Inequality (19) follows from (23). Inequality (19) implies (20) since, for $1 \leq i \leq n-2$, $\sqrt{\frac{j-1}{i+1}}$ takes its minimum when $i=n-2$. This completes the proof of the lemma.

We have found all critical points for the max-min problem. Namely, for each $i$, $1 \leq i \leq n-1$, there are $\binom{n}{i}$ critical points of type $i$; the corresponding critical value is given by (16). (A critical point of type i can be described geometrically as follows. Consider a ray from the centroid $C_0$ of the (n-1)-simplex S given in (2) to the centroid of a k-dimensional face of S ($0 \leq k \leq n-2$). The intersection of this ray with the sphere of radius $\sigma\sqrt{n}$ centered at $C_0$ is a critical point of type k+1.) In addition, there are the critical points **x** where **grad** $G(\mathbf{x})=\mathbf{0}$ (i.e. points **x** which have two or more of their coordinates equal to 0), and there is the critical point $(\mu,...,\mu)$ when $\sigma = 0$.

Now consider part (a) of Theorem 1 where $0 \leq \sigma/\mu < 1/\sqrt{n-1}$. We may assume that $0 < \sigma/\mu < 1/\sqrt{n-1}$ since, as remarked earlier, if $\sigma = 0$ then (3) is trivial. Consider the set of points **x** in n-space whose coordinates have the given mean $\mu$ and variance $\sigma^2$. By Lemma 1(a), these **x** have all their coordinates positive. By compactness, the function $G(\mathbf{x})$ restricted to this set attains a maximum and minimum. Therefore the maximum and minimum must occur among the critical values $G(\mathbf{x_i})$ given by (16). Since $0 < t = \sigma/\mu < 1/\sqrt{n-1}$ we see from (17) and (20) that for all **x** with the given mean and variance

(27) $$P_{n-1}(t) = \frac{G(\mathbf{x}_{n-1})}{\mu^n} \leq \frac{G(\mathbf{x})}{\mu^n} \leq \frac{G(\mathbf{x}_1)}{\mu^n} = P_1(t)$$

which proves (3).

Now consider part (b) of Theorem 1. Here we are given an n-vector $(x_1, x_2, \ldots, x_n)$ where the $x_i$ are positive with mean $\mu$ and variance $\sigma^2$. By Lemma 1(b), $0 \leq \sigma/\mu < \sqrt{n-1}$. As before, we can dispense with the trivial case $\sigma = 0$. We want to find the maximum of $G(\mathbf{y})$ among all n-vectors **y** whose coordinates are positive and have the given $\mu$ and $\sigma^2$. By compactness G attains a maximum on the intersection

(28) $$\{A(\mathbf{y})= n\mu\} \cap \{V(\mathbf{y})=n\sigma^2\} \cap \{y_1 \geq 0, y_2 \geq 0,\ldots,y_n \geq 0\}.$$

This maximum must occur at a point $\mathbf{y_0}$ with positive coordinates. Therefore this maximum is a local maximum for the max-min problem (7) and hence occurs at a critical point. That is, $\mathbf{y_0} = \mathbf{x_i}$ for some $i$, where $\mathbf{x_i}$ is a critical point of type *i*. Although it is possible in this case (b) for $G(\mathbf{x_1}) < G(\mathbf{x_i})$ for some i, we can make use of the observation

(29)   $G(\mathbf{x_1}) = \max_i \{G(\mathbf{x_i}): \mathbf{x_i}$ is a critical point having all coordinates positive$\}$

which follows from (15), (16), and (19). Therefore $\mathbf{y_0} = \mathbf{x_1}$ for some critical point $\mathbf{x_1}$ of type 1. Hence $x_1 x_2 \cdots x_n \leqq G(\mathbf{x_1}) = \mu^n P_1(\sigma/\mu)$, which establishes the upper bound in (3) in case (b) of Theorem 1. This completes the proof of Theorem 1.

## Section 2

As explained in the papers of Aldaz [1,2], Becker [4], Estrada [6], and Markowitz [9], generalizations of the arithmetic and geometric means inequality which involve only the variance of the sequence are the most useful in applications to economics and finance. The paper Becker [4] contains a discussion, with historical references, of the heuristics behind the approximation $R_A - R_G \approx \sigma^2/2$, where $R_A$ and $R_G$ denote the arithmetic and geometric mean. Markowitz [9] considers five different mean-variance approximations for the geometric mean and compares their accuracy for sequences of historical economic data. Theorem 2.4 of Aldaz [2] contains a general inequality involving weighted means and generalized variances which is optimal within its class. When this general inequality is specialized by setting the weights $\alpha=(1/n,1/n,...,1/n)$ and s=2 the result becomes, for nonnegative sequences,

(30) $$R_A - R_G \leq n\sigma.$$

Theorem 1 can be applied to obtain a similar form of upper bound for $R_A - R_G$. Indeed, the lower bound in (3) implies $(\mu - \sigma\sqrt{n-1})^n \leq x_1 x_2 \cdots x_n$, and hence $R_A - R_G \leq \sqrt{n-1}\ \sigma$. We record this result as a corollary to Theorem 1:

**Corollary 1**. *Fix* $n \geqq 2$. *If* $x_1, x_2, \ldots, x_n$ *is a positive sequence with mean* $\mu$ *and variance* $\sigma^2$ *then*
(31) $$\mu - (x_1 x_2 \cdots x_n)^{1/n} \leq \sqrt{n-1}\ \sigma.$$

Aldaz has shown there can be no similar lower bound; i.e., there does not exist a constant k>0 such that $k\sigma \leq \mu - (x_1 x_2 \cdots x_n)^{1/n}$ is valid for all positive sequences $x_1, x_2, \ldots, x_n$ with mean $\mu$ and variance $\sigma^2$ (see Example 2.1 of [1]).

A number of generalizations of the AM-GM inequality in the literature involve properties other than variance. Cartwright and Field [5] prove an inequality involving weighted arithmetic means, variance, and upper and lower bounds for the sequence. In the special case of equal weights their result reduces to

$$\text{(32)} \qquad \frac{\sigma^2}{2b} \leq R_A - R_G \leq \frac{\sigma^2}{2a}$$

where a,b denote lower and upper bounds, respectively, for the positive sequences being considered (see also Wang[14]). For easier comparision with (3) we can rewrite the inequality (32) of Cartwright and Field in the equivalent form

$$\text{(33)} \qquad (\mu - \frac{\sigma^2}{2a})^n \leq x_1 x_2 \cdots x_n \leq (\mu - \frac{\sigma^2}{2b})^n.$$

Alzer [3] proves a refinement of the inequality of Cartwright and Field [5] which incorporates variance but also retains the bounds a,b defined above. Tung [12] derives inequalities depending only on the bounds a,b—his inequalities do not involve variance. Meyer [10] extends those results to the harmonic mean. Aldaz [1] makes use of the variance of the square roots of the terms of the sequence, and in [2] he extends those results to more general weights and variances. Loewner and Mann [8] derive an upper bound which involves the maximum and minimum of $x_i/\mu$ and does not incorporate variance.

**Remark 2.** We illustrate the relevance of Theorem 1 to finance. Consider an investment in a certain asset A. Suppose that for n consecutive time periods the investment returns are $r_1, r_2, \ldots, r_n$ (-1< $r_i$). For example, if the time period is years and if asset A returned 6% in the i-th year then $r_i$ = .06; if it lost 6% that year then $r_i$ = - .06.

An initial investment of $1 in asset A will be worth $X_n$ at the end of the n-th year, where $X_n = (1+r_1) \cdots (1+r_n)$. Suppose the sequence $r_1, r_2, \ldots, r_n$ has mean $\mu_n$ and variance $\sigma_n^2$. Then the sequence $1+r_1, \ldots, 1+r_2$ will have mean $1+\mu_n$, and variance $\sigma_n^2$. The terms of this sequence are positive since -1< $r_i$. By Theorem 1

$$\text{(34)} \quad \left(1+\mu_n - \sigma_n\sqrt{n-1}\right)\left(1+\mu_n + \frac{\sigma_n}{\sqrt{n-1}}\right)^{n-1} \leq X_n \leq \left(1+\mu_n + \sigma_n\sqrt{n-1}\right)\left(1+\mu_n - \frac{\sigma_n}{\sqrt{n-1}}\right)^{n-1}.$$

Remark 1 shows—perhaps unexpectedly—that for fixed $\mu_n$, $\sigma_n$ with $\sigma_n > 0$, the best investment outcome $X_n$ occurs when all returns but one are identical and below the mean $\mu_n$; the worst outcome X occurs when all returns but one are identical and above the mean $\mu_n$.

Before making the initial investment an investor can estimate the mean and variance of returns for the asset A from its historical performance record; let $\mu_0$ and $\sigma_0$ be the values so obtained. For example, suppose asset A is the S&P 500 index and suppose the unit of time is days. Based on the historical record from 3 January 1950 through 31 July 2012 it has been estimated that the daily returns on this asset will have a mean $\mu_0$ of 1.0003 and a standard deviation $\sigma_0$ of .0098 (cf. [13]).

Suppose one expects that $\mu_n$, $\sigma_n$ will be close to their estimated values $\mu_0$, $\sigma_0$; say $|\mu_n - \mu_0| < \varepsilon$ and $|\sigma_n - \sigma_0| < \varepsilon$, $(\varepsilon > 0)$. Then Equation (34) will provide the following estimate for $X_n/(1+\mu_n)^n$, the ratio of outcomes for an n-term investment in a risk free asset with the same mean:

$$(35) \quad \frac{X_n}{(1+\mu_n)^n} \leq \left(1 + \frac{(\sigma_0+\epsilon)\sqrt{n-1}}{(1+\mu_0-\epsilon)}\right)\left(1 - \frac{(\sigma_0-\epsilon)}{(1+\mu_0+\epsilon)\sqrt{n-1}}\right)^{n-1}.$$

Note that if $\sigma_0 - \varepsilon > 0$ and $1 + \mu_0 - \varepsilon > 0$ then the right hand side of Equation (35) will tend to 0 and n tends to infinity.

Department of Mathematics
University of California, San Diego
brodin@ucsd.edu________________________________